\documentstyle[12pt]{article}

\newcommand{\ptl}{\partial}
\newcommand{\dps}{\displaystyle}

\newtheorem{theorem}{\bf Theorem}[section]
\newtheorem{remark}{\bf Remark}[section]

\begin{document}

\begin{center}
{\Large \bf Gap of the First Two Eigenvalues of the Schr\"odinger 

\medskip

Operator with Nonconvex Potential}

\vskip 1cm

{\large \bf Shing-Tung Yau} 

\vskip .75cm

{\it Dedicated to Manfredo do Carmo on his 80th Birthday}

\end{center}

\vskip 1cm

In this essay, I will extend my previous work \cite{Yau} on operators
whose potential is nonconvex. In particular, the results given here can be
applied to the double well potential. I define an invariant associated
to the potential in \S 4.  It defines a distance between the point
where $\textstyle \frac{u_2}{u_1}$ achieves its maximum, $\sup\, \frac{u_2}{u_1}\,$, to the point where $\frac{u_2}{u_1}=\varepsilon \sup\, \frac{u_2}{u_1}$. Here $u_i$ are the eigenfunctions of the Schr\"odinger operator. I will show how the gap
$\lambda_2-\lambda_1$ can be estimated from below in terms of this
distance.  Theorem 6.1 is the main theorem of this essay.  It is a very interesting problem to locate the maximum of $\frac{u_2}{u_1}$ and its zeroes.  The upper bound of $\lambda_2-\lambda_1$ depends on the choice of a good trial function, and I shall come back to this question in the future.

This line of research on gradient estimates started from my work on
bounded harmonic functions \cite{Yau75} and the method was used by Peter Li \cite{Li} and Li-Yau \cite{LY} for the Laplacian of a manifold. Li-Yau \cite{LY} also applied it to the Schr\"odinger operator where a distance function similar to the one used here was
introduced.

The Li-Yau type distance function was also used by Perelman in his
famous work \cite{Perelman}.

In the Li-Yau's approach of estimating the first eigenvalue of the
Laplacian, it was conjectured by Li-Yau and proved by Zhong-Yang \cite{Zhong Yang}
that if $d$ is the diameter of a manifold with nonnegative Ricci
curvature, then $\lambda_1 \, d^2$ has an universal lower bound
which is achieved when the manifold is a circle.

Convex domain and convex potential can be considered as an analogue
of manifold with nonnegative curvature. In Singer-Wong-Yau-Yau \cite{SWYY}, we improved on the log concavity result of Brascamp-Lieb \cite{BL} and proved that $(\lambda_2-\lambda_1) d^2$ has a universal lower bound.  It is natural for us to expect that the interval will give this optimal lower estimate.

I would like to dedicate this work to my friend Manfredo do Carmo whose
works on minimal surfaces are very original and influential.

\vskip 1cm

\section{Generalized log concavity of the first eigenfunction}

In \cite{SWYY,Yau}, I used method of continuity to generalize the log
concavity result of Brascamp and Lieb \cite{BL} when the potential
is convex. I generalize it further in this section.

Let $u_1$ be the first eigenfunction of the operator $-\Delta +V$
on a domain $\Omega_1$ with zero boundary valued data. Let
$\varphi=-\log u_1$. Then we have the following theorem:

\begin{theorem}
The Hessian of $\varphi$ has eigenvalue greater than $g(x)$ where
\begin{eqnarray*}
\begin{array}{rcl}
g(x) = & \sup \{f(x)| & f \ {\rm is\ a\ bounded\ smooth\ function\
defined\ on}\ \Omega \ {\rm so\ that\ the}\\
& & {\rm lowest\ eigenvalue\ of\ the\ Hessian\ of}\ V \ {\rm
plus}\ \Delta f \ {\rm is\ greater\ than}\ f^2 \}
\end{array}
\end{eqnarray*}
\end{theorem}

\vskip 5mm

\noindent {\bf Proof} \ Differenting the equation
\begin{eqnarray}
\Delta\varphi = |\nabla\varphi|^2 - V+\lambda_1
\end{eqnarray}
we obtain
\begin{eqnarray}
\Delta(\varphi_{ii} -f)=\sum\varphi^2_{ji} - V_{ii} - \Delta f
\end{eqnarray}
where $\dps \varphi_{ji}=\frac{\ptl^2 u}{\ptl x_j \, \ptl x_i}$
and $\dps V_{ii}=\frac{\ptl^2V}{\ptl x^2_i}$.

Minimizing $\varphi_{ii}-f$, we obtain
\begin{eqnarray}
\varphi_{ji}&=&0 \qquad {\rm for} \qquad j\not= i \\
\Delta(\varphi_{ii}-f) &\geq & 0
\end{eqnarray}
Hence at such a point,
\begin{eqnarray}
\varphi^2_{ii} \geq V_{ii} + \Delta f
\end{eqnarray}

Note that the continuity argument was introduced by me and discussed in \cite{SWYY}.  (I had lectured on this in 1979 as was noted in \cite{Korevaar}.)  It can be applied
in the following way.

Replace the potential $V$ by $V_t = ||x||^2 + tV\,$.  When $t=0$, the
theorem is obviously true. We have to prove the theorem for all
$t>0$.

Suppose the theorem is true for $t<t_0$. Then at $t=t_0$,
$\varphi_{ii} \geq f$ and there is a function $f$, depending on
$t$, so that $(V_t)_{ii}+ \Delta f> f^2$ and at some point,
$\varphi_{ii}=f$. By (1.5), we obtain
\begin{eqnarray}
f^2 \geq (V_t)_{ii}+ \Delta f
\end{eqnarray}
This contradict the choice of $f$.

Hence $t$ can be arbitrarily large and we conclude that Theorem 1.1
holds.

\begin{remark}
The argument of Theorem 1.1 can be generalized to manifolds with
negative curvature.
\end{remark}

\vskip 1cm

\section{Gradient estimate of first eigenfunction}

\setcounter{equation}{0}

Let $\Omega$ be a convex subdominant. Then we can choose a smooth
nonnegative function $\rho$ with compact support.

Let
\begin{eqnarray}
G=\rho^2 (V+\alpha)^{-1} \ \Delta\varphi
\end{eqnarray}
where $\alpha$ is a constant to be chosen later. Then according to
(1.1), we find
\begin{eqnarray}
\nabla G=G(2\nabla\log\rho - \nabla\log (V+\alpha)) + \rho^2
(V+\alpha)^{-1} \ (2\nabla\varphi \cdot \nabla\nabla\varphi -
\nabla V)
\end{eqnarray}
\begin{eqnarray}
\begin{array}{rcl}
\Delta G & = & \nabla G (2\nabla\log\rho - \nabla\log (V+\alpha))
+ G (2\Delta\log\rho - \Delta\log (V+\alpha))\\
& & + 2\nabla G \cdot \nabla\log\rho - 4 |\nabla\log\rho|^2 \ G -
\nabla G \cdot \nabla\log (V+\alpha)\\
& & +4G\nabla\log (V+\alpha) \cdot \nabla\log\rho - G|\nabla\log
(V+\alpha)|^2\\
& & + \rho^2 \ (V+\alpha)^{-1} \ [2|\nabla\nabla\varphi|^2 +
4\nabla\varphi \cdot \nabla\nabla\varphi \cdot \nabla\varphi -
2\nabla\varphi \cdot \nabla V - \Delta V]\\
& = & \nabla G (4\nabla\log\rho - 2\nabla\log (V+\alpha) +
2\nabla\varphi) \\
& & + G ( 2\Delta\log\rho - \Delta\log (V+\alpha) - 4
|\nabla\log\rho|^2 - |\nabla\log (V+\alpha)|^2 \\
& & \hskip 1cm - 4\nabla\log\rho \cdot \nabla\varphi + 2\nabla\log
(V+\alpha) \cdot
\nabla\log\rho )\\
& & + \rho^2 (V+\alpha)^{-\frac{1}{2}} \ (|\nabla\nabla\varphi|^2
- \Delta V)
\end{array}
\end{eqnarray}
By (1.1)
\begin{eqnarray}
\begin{array}{rcl}
|\nabla\nabla\varphi|^2 & \geq & \dps \frac{1}{n}
|\Delta\varphi|^2\\[3mm]
& = & \dps \frac{1}{n} (V+\alpha)^2 \ \rho^4 \ G^2
\end{array}
\end{eqnarray}
Hence at the point where $G$ achieves its maximum,
\begin{eqnarray}
\begin{array}{rcl}
0 & \geq & \dps \frac{2}{n} G^2 - \rho^4 (V+\alpha)^{-2} \ \Delta
V \\
& & + \rho^2 \ (V+\alpha)^{-1} \ G [ 2\Delta\log\rho - \Delta\log
(V+\alpha) - 4|\nabla\log\rho|^2 \\
& & \hskip 2cm - |\nabla\log (V+\alpha)|^2 -4\log\rho \cdot
\nabla\varphi + 2\nabla\log (V+\alpha) \cdot \nabla\log\rho]
\end{array}
\end{eqnarray}

Note that
\begin{eqnarray}
\begin{array}{rcl}
\rho^2 (V+\alpha)^{-1} \ |\nabla\varphi|^2 & \leq & \rho^2 \
(V+\alpha)^{-1} \ (|\nabla\varphi|^2 - V - \lambda_1) + \rho^2 \
(V+\alpha)^{-1} \ (V-\lambda_1)\\
& = & G+\rho^2 \ (V+\alpha)^{-1} \ (V-\lambda_1)
\end{array}
\end{eqnarray}
Hence
\begin{eqnarray}
\begin{array}{rcl}
\dps \frac{2}{n} G^2 & \leq & \rho^4 \ (V+\alpha)^{-2} \ \Delta
V\\
& & -\rho^2 \ (V+\alpha)^{-1} \ G [2\Delta\log\rho - \Delta\log
(V+\alpha) -4|\nabla\log\rho|^2 \\
& & \hskip 1cm - |\nabla\log (V+\alpha)|^2 + 2\nabla\log
(V+\alpha) \cdot \nabla\log\rho] \\
& & +4\rho (V+\alpha)^{-\frac{1}{2}}
\end{array}
\end{eqnarray}

Therefore,

\noindent either
\begin{eqnarray}
G\leq \frac{4}{n^2} \ \rho^2 (V+\alpha)^{-1} \ |\nabla\log\rho|^2
\end{eqnarray}
or
\begin{eqnarray}
G \leq 3\sqrt{\frac{n}{2}} \ \rho^2 \ (V+\alpha)^{-1} \
\sqrt{(\Delta V)_+}
\end{eqnarray}
or
\begin{eqnarray}
\begin{array}{rcl}
G & \leq & \dps \frac{3n}{2} \ \rho^2 \ (V+\alpha)^{-1} \
[-2\Delta\log\rho + 4|\nabla\log\rho|^2 + \Delta\log (V+\alpha)\\
& & \hskip 1cm + |\nabla\log (V+\alpha)|^2 - 2\nabla\log
(V+\alpha) \cdot \nabla\log\rho + 4 (V-\lambda_1)^{\frac{1}{2}}]
\end{array}
\end{eqnarray}

\vskip 5mm

\begin{theorem}
Let $u_1$ be a positive solution of $(\Delta+V)u_1=\lambda_1\,u$
and $\varphi=-\log u_1$. Then
\begin{eqnarray}
\begin{array}{rcl}
\rho^2 \ (V+\alpha)^{-1}&&\!\!\!\!\!\!\!\! \!\!\!\!\!(|\nabla\varphi|^2- V+\lambda_1) \\
& \leq&  10n \ \rho^2 \ (V+\alpha)^{-1} \ (|\nabla\log\rho|^2 +
|\Delta\log\rho|)\\
& & \dps + \frac{3n}{2} \ \rho^2 \ (V+\alpha)^{-1} \ \left\{
(\Delta\log (V+\alpha))_+ + 2|\nabla\log (V+\alpha)|^2 +
4(V-\lambda_1)^{\frac{1}{2}} \right\}
\end{array}
\end{eqnarray}
In particular,
\begin{eqnarray*}
\begin{array}{rcl}
\rho^2 \ (V+\alpha)^{-1}|\nabla\varphi|^2
 &\leq & \dps \frac{\sup V-\lambda_1}{\sup V+\alpha} + 10n^2 \sup
(V+\alpha)^{-1} \ \rho^2 \ (|\nabla\log\rho|^2 +
|\Delta\log\rho|)\\
& & + 3n\sup \left[ (V+\alpha)^{-2} \ ((\Delta V)_+) +
(V+\alpha)^{-3} \ |\nabla V|^2 \right]\\
& & \dps + 6n\sup (V+\alpha)^{-\frac{1}{2}} \ \left( \frac{\sup
V-\lambda_1}{\sup V+\alpha} \right)^{\frac{1}{2}}
\end{array}
\end{eqnarray*}
\end{theorem}

\vskip 1cm

\section{Gradient estimate for $\dps \frac{u_2}{u_1}$}

\setcounter{equation}{0}

Let $u_2$ be the second eigenfunction of $-\Delta+V$ on $\Omega$.
Let $\dps u=\frac{u_2}{u_1}$.

\noindent Then
\begin{eqnarray}
\Delta u=-(\lambda_2-\lambda_1)u +2\nabla\varphi \, \nabla u
\end{eqnarray}

Let $\dps c>\sup_{\Omega_2} \frac{u_2}{u_1}$ and $\psi=-\ln(c-u)$.

\noindent Then
\begin{eqnarray}
\Delta\psi=(\lambda_2-\lambda_1)(1-ce^\psi)+2\nabla\varphi
\nabla\psi+|\nabla\psi|^2.
\end{eqnarray}

Let
\begin{eqnarray}
F=\rho^2 (V+\alpha)^{-1} \left[ |\nabla\psi|^2 +
(\lambda_2-\lambda_1) (1-ce^\psi) \right]
\end{eqnarray}
when $\alpha>0$ is a constant to be chosen.

\noindent Then
\begin{eqnarray}
\nabla F = 2F \nabla\log\rho - F\nabla\log (V+\alpha) + \rho^2
(V+\alpha)^{-1} \left[ 2\nabla\psi \, \nabla\nabla\psi -
c(\lambda_2-\lambda_1) e^\psi \nabla\psi \right]
\end{eqnarray}
and
\begin{eqnarray}
\begin{array}{rcl}
\Delta F & = & \nabla F \left( 2\nabla\log\rho - \nabla\log
(V+\alpha) \right) + F \left( 2\Delta\log\rho - \Delta\log
(V+\alpha) \right)\\
& & + 4\rho (V+\alpha)^{-1} \nabla\rho \cdot \nabla\psi \cdot
\nabla\nabla\psi - 2\rho^2 (V+\alpha)^{-2} \ \nabla V \cdot
\nabla\nabla\psi \cdot \nabla\psi \\
& & + 2\rho^2 (V+\alpha)^{-1} \ |\nabla\nabla\psi|^2 + 2\rho^2
(V+\alpha)^{-1} \ \nabla\psi \cdot \nabla(\Delta\psi) \\
& & - 2c\rho (V+\alpha) (\lambda_2-\lambda_1) e^\psi \ \nabla\rho
\cdot \nabla\psi + c\rho^2 (V+\alpha)^{-2} (\lambda_2-\lambda_1)
e^\psi \ \nabla V \cdot \nabla\psi \\
& & - c\rho^2 (V+\alpha)^{-1} (\lambda_2-\lambda_1) e^\psi \
|\nabla\psi|^2 \\
& & -c\rho^2 (V+\alpha)^{-1} (\lambda_2-\lambda_1) e^\psi \left[
|\nabla\psi|^2 + 2\nabla\varphi \cdot \nabla\psi +
(\lambda_2-\lambda_1) (1-ce^\psi) \right] \\
& = & \nabla F(2\nabla\log\rho - \nabla\log (V+\alpha)) +
F(2\Delta\log\rho - \Delta\log (V+\alpha))\\
& & + 2\nabla F \cdot \nabla\log\rho - 4F |\nabla\log\rho|^2 + 2F
\nabla\log (V+\alpha) \cdot \nabla\log\rho \\
& & + 2c\rho^2 (V+\alpha)^{-1} (\lambda_2-\lambda_1) e^\psi \
\nabla\psi \cdot \nabla\log\rho - (V+\alpha)^{-1} \ \nabla F \cdot
\nabla V \\
& & + 2F\nabla\log\rho \cdot \nabla\log (V+\alpha) - F|\nabla\log
(V+\alpha)|^2 \\
& & -c\rho^2 (V+\alpha)^{-2} (\lambda_2-\lambda_1) e^\psi \
\nabla\psi \cdot \nabla V + 2\rho^2 (V+\alpha)^{-1}
|\nabla\nabla\psi|^2 \\
& & +2\rho^2 (V+\alpha)^{-1} (2\nabla\psi \cdot \nabla\nabla\psi
\cdot \nabla\psi + 2\nabla\psi \cdot \nabla\nabla\psi \cdot
\nabla\psi \\
& & \hskip 3cm + 2\nabla\psi \cdot \nabla\psi \cdot
\nabla\nabla\psi -c(\lambda_2-\lambda_1)e^\psi |\nabla\psi|^2) \\
& & -2c\rho (V+\alpha) (\lambda_2-\lambda_1) e^\psi \nabla\psi
\cdot \nabla\psi + c\rho^2 (V+\alpha)^{-2} (\lambda_2-\lambda_1)
e^\psi \ \nabla V \cdot \nabla\psi \\
& & -c\rho^2 (V+\alpha)^{-1} (\lambda_2-\lambda_1) e^\psi
|\nabla\psi|^2\\
& & -c\rho^2 (V+\alpha)^{-1} (\lambda_2-\lambda_1) e^\psi \left[
|\nabla\psi|^2 + 2\nabla\varphi \cdot \nabla\psi +
(\lambda_2-\lambda_1) (1-ce^\psi) \right] \\
& = & \nabla F(4\nabla\log\rho - 2\nabla\log (V+\alpha)) + F (
2\Delta\log\rho - \Delta\log (V+\alpha) - 4 |\nabla\log\rho|^2 \\
& & \hskip 3cm + 4\nabla\log (V+\alpha) \cdot \nabla\log\rho -
|\nabla\log (V+\alpha)|^2 ) \\
& & +2\rho^2 (V+\alpha)^{-1} \ |\nabla\nabla\psi|^2 + 2 \nabla F
\cdot \nabla\psi - 2F \nabla\log\rho \cdot \nabla\psi \\
& & -2F\nabla\log (V+\alpha) \cdot \nabla\psi + 4\rho^2
(V+\alpha)^{-1} \ \nabla\psi \cdot \nabla\nabla\psi \cdot
\nabla\psi + 2\nabla F \cdot \nabla\varphi \\
& & - 2F \nabla\log\rho \cdot \nabla\varphi + 2F\nabla\log
(V+\alpha) \nabla\varphi + 2c\rho^2 (V+\alpha)^{-1}
(\lambda_2-\lambda_1) e^\psi \ \nabla\psi \cdot \nabla\varphi \\
& & -c\rho^2 (V+\alpha)^{-1} (\lambda_2-\lambda_1) e^\psi
|\nabla\psi|^2 \\
& & -c\rho^2 (V+\alpha)^{-1} (\lambda_2-\lambda_1) e^\psi \left[
|\nabla\psi|^2 + 2\nabla\varphi \cdot \nabla\psi +
(\lambda_2-\lambda_1) (1-ce^\psi) \right] \\
& = & \nabla F (4\nabla\log\rho - 2\nabla\log (V+\alpha) +
2\nabla\psi + 2\nabla\psi) \\
& & + F(2\Delta\log\rho - \Delta\log (V+\alpha) -
4|\nabla\log\rho|^2 \\
& & \hskip 1cm + 4\nabla\log (V+\alpha) \cdot \nabla\log\rho -
|\nabla\log (V+\alpha)|^2 - 2\nabla\log\rho \cdot \nabla\psi \\
& & \hskip 1cm - 2 \nabla\log (V+\alpha) \cdot \nabla\psi - 2
\nabla\log\rho \cdot \nabla\varphi + 2\nabla\log (V+\alpha) \cdot
\nabla\varphi) \\
& & + 2\rho^2 (V+\alpha)^{-1} \ |\nabla\nabla\psi|^2 + 4\rho^2
(V+\alpha)^{-1} \ \nabla\psi \cdot \nabla\nabla\psi \cdot
\nabla\psi \\
& & -2c\rho^2 (V+\alpha)^{-1} (\lambda_2-\lambda_1) e^\psi \
|\nabla\psi|^2 - c(\lambda_2-\lambda_1)^2 \ \rho^2 \
(V+\alpha)^{-1} \ e^\psi (1-ce^\psi)
\end{array}
\end{eqnarray}

Now
\begin{eqnarray}
|\nabla\nabla\psi|^2 \geq \frac{1}{n} (\Delta\psi)^2
\end{eqnarray}
and
\begin{eqnarray}
\Delta\psi - 2\nabla\psi \cdot \nabla\varphi = \rho^{-2} \
(V+\alpha) F
\end{eqnarray}

\noindent Hence
\begin{eqnarray}
\begin{array}{rcl}
& & \rho^2 (V+\alpha)^{-1} \ |\nabla\nabla\psi|^2 \\
& \geq & \rho^2 (V+\alpha)^{-1} \left[ (V+\alpha)^2 \ \rho^{-4} \
F^2 + 4\rho^{-2} \ (V+\alpha) F\nabla\psi \cdot \nabla\varphi + 4
(\nabla\psi \cdot \nabla\varphi)^2 \right]
\end{array}
\end{eqnarray}

When $F$ achieves its maximum, $\nabla F=0$ and $\Delta F \leq 0$.
Therefore,
\begin{eqnarray}
\begin{array}{rcl}
0 & \geq & F(2\Delta\log\rho - \Delta\log (V+\alpha) - 4
|\nabla\log\rho|^2 + 4\nabla\log (V+\alpha) \cdot \nabla\log\rho
\\
& & \hskip 1cm - |\nabla\log (V+\alpha)|^2 - 2\nabla\log\rho \cdot
\nabla\psi - 2\nabla\log (V+\alpha) \cdot \nabla\psi \\
& & \hskip 1cm - 2 \nabla\log\rho \cdot \nabla\varphi +
2\nabla\log (V+\alpha) \cdot \nabla\varphi) \\
& & + \dps \frac{2}{n} \rho^{-2} \ (V+\alpha) F^2 + \frac{8}{n} F
\ \nabla\psi \cdot \nabla\varphi + \frac{8}{n} \rho^2
(V+\alpha)^{-1} \ (\nabla\psi \cdot \nabla\varphi)^2 \\
& & +4\rho^2 (V+\alpha)^{-1} \ \nabla\psi \cdot
\nabla\nabla\varphi \cdot \nabla\psi \\
& & -2c\rho^2 (V+\alpha)^{-1} \ (\lambda_2-\lambda_1) e^\psi
\left[ |\nabla\psi|^2 + (\lambda_2-\lambda_1) (1-ce^\psi) \right]
\\
& & +c(\lambda_2-\lambda_1)^2 \ \rho^2 \ (V+\alpha)^{-1} \ e^\psi
(1-ce^\psi)
\end{array}
\end{eqnarray}

Note that
\begin{eqnarray}
\rho^2 (V+\alpha)^{-1} \ \nabla\psi \cdot \nabla\nabla\varphi
\cdot \nabla\psi \geq \inf \left( \rho^2 (V+\alpha)^{-1} \ g(x)
\right) F
\end{eqnarray}
where $g$ is defined in Theorem 1.1.

\noindent Hence
\begin{eqnarray}
\begin{array}{rcl}
0 & \geq & \dps \frac{2}{n} F^2 + 4\inf \left(\rho^2 \
(V+\alpha)^{-1} \ g(x) \right) F\\
& & \dps - F \left( \rho^2 \, (V+\alpha)^{-1} \, |\nabla\psi|^2
\right) \left[ \frac{8}{n} \left( \rho^2 \ (V+\alpha)^{-1} \
|\nabla\varphi|^2 \right)^{\frac{1}{2}} + 2 \left( \rho^2 \
(V+\alpha)^{-1} \
|\nabla\log\rho|^2 \right) \right.\\
& & \dps \hskip 1cm + 2 \left( \rho^2 \ (V+\alpha)^{-1} \
|\nabla\log
(V+\alpha)|^2 \right)^{\frac{1}{2}} \bigg]\\
& & \dps +\rho^2 \, (V+\alpha)^{-1} \, F \left[ 2\Delta\log\rho -
7|\nabla\log\rho|^2 - 2|\nabla\varphi|^2 - \Delta\log (V+\alpha)- |\nabla\log (V+\alpha)|^2 \right] \\
& & \dps - 2c\rho^2 (V+\alpha)^{-1} \ (\lambda_2-\lambda_1) e^\psi
\, F + c(\lambda_2-\lambda_1)^2 \, \rho^\psi \ (V+\alpha)^{-2} \
e^\psi \ (1-ce^\psi)
\end{array}
\end{eqnarray}

By definition of $F$, either
\begin{eqnarray}
\rho^2 \ (V+\alpha)^{-1} \ |\nabla\psi|^2 \leq 2F
\end{eqnarray}
or
\begin{eqnarray}
F \leq (\lambda_2-\lambda_1) (V+\alpha)^{-1} \ (1-ce^\psi)
\end{eqnarray}

Let us assume (3.12) first. In that case, either
\begin{eqnarray}
\begin{array}{rcl}
F^{\frac{1}{2}} & \leq & \dps \frac{3n}{2} \left[ \frac{8}{n}
\left( \rho^2 \ (V+\alpha)^{-1} \ |\nabla\varphi|^2
\right)^{\frac{1}{2}} + 2 \left( \rho^2 \ (V+\alpha)^{-1} \
|\nabla\log\rho|^2
\right)^{\frac{1}{2}} \right] \\
& & \dps + 2 \left( \rho^2 \ (V+\alpha)^{-1} \ |\nabla\log
(V+\alpha)|^2 \right)^{\frac{1}{2}}
\end{array}
\end{eqnarray}
or
\begin{eqnarray}
\begin{array}{rcl}
& & \dps F+6n \inf \left( \rho^2 \ (V+\alpha)^{-1} \ g(x)
\right)\\
& \leq & \dps \frac{3n}{2} \rho^2 \ (V+\alpha)^{-1} \bigg[
2\Delta\log\rho + 7 |\nabla\log\rho|^2 + \Delta\log (V+\alpha) \\
& & \hskip 1cm \dps + |\nabla\log (V+\alpha)|^2 +
2|\nabla\varphi|^2 + 2c (\lambda_2-\lambda_1) e^\psi \bigg]
\end{array}
\end{eqnarray}
or
\begin{eqnarray}
F^2 \leq \frac{3n}{2} c (\lambda_2-\lambda_1)^2 \ \rho^2 \
(V+\alpha)^{-2} \ e^\psi \ |(1-ce^\psi)|
\end{eqnarray}

From this, we conclude:
\begin{theorem}
Let $\Omega$ be a domain and $\rho$ be smooth function with
compact support in $\Omega$. Let $u_i$ be smooth function
satisfying the equation $(-\Delta+V) u_i = \lambda_i \, u_i$ so
that $u_1 >0$ and $\varphi=-\log u_1$ satisfies the conclusion of
Theorem 1.1. Let $c$ be a constant so that $c>\sup u$ where $\dps
u=\frac{u_2}{u_1}$. Let $\psi=-\log (c-u)$. Then one of the
following inequalities hold:
\begin{enumerate}
\item[(1)]
\begin{eqnarray}
\rho^2 (V+\alpha)^{-1} \left[ |\nabla\psi|^2 +
(\lambda_2-\lambda_1) (1-ce^\psi) \right] \leq
(\lambda_2-\lambda_1) \sup (V+\alpha)^{-1} (1-ce^\psi) \quad
\end{eqnarray}

\item[(2)]
\begin{eqnarray}
\begin{array}{rcl}
& & \dps \rho^2 \ (V+\alpha)^{-1} \left[ |\nabla\psi|^2 +
(\lambda_2-\lambda_1) (1-ce^\psi) \right] \\
& \leq & \dps 144 \left( \frac{\sup V-\lambda_1}{\sup
V+\alpha}\right) + 20n^2 \, \sup (V+\alpha)^{-1} \ \rho^2 \ \left(
|\nabla\log\rho|^2 + |\Delta\log\rho| \right) \\
& & \dps + 4n \, \sup \left( (V+\alpha)^{-2} \ (\Delta V)_+ +
(V+\alpha)^{-3} \ |\nabla V|^2 \right) \\
& & \dps +6n \, \sup (V+\alpha)^{-\frac{1}{2}} \ \left( \frac{\sup
V-\lambda_1}{\sup V+\alpha} \right)^{\frac{1}{2}}
\end{array}
\end{eqnarray}

\item[(3)]
\begin{eqnarray}
\begin{array}{cl}
& \dps \rho^2 \left[ (V+\alpha)^{-1} \ |\nabla\psi|^2 +
(\lambda_2-\lambda_1) (1-ce^\psi) \right] + 6n \inf \left( \rho^2
\, (V+\alpha)^{-1} \, g(x) \right)\\
\leq & \dps 20n^2 \, \sup (V+\alpha)^{-1} \ \left( |\nabla\rho|^2
+ \rho|\Delta\rho| \right)\\
& \dps +10n \, \sup \rho^2 \left[ (V+\alpha)^{-2} \, (\Delta V)_+
+ (V+\alpha)^{-3} \, |\nabla V|^2 \right]\\
& \dps +10 \left( \frac{\sup V-\lambda_1}{\sup V+\alpha}\right) +
6n \, \sup (V+\alpha)^{-\frac{1}{2}} \, \left( \frac{\sup
V-\lambda_1}{\sup V+\alpha}\right)^{\frac{1}{2}}\\
& \dps +3nc (\lambda_2-\lambda_1) \, \sup \rho^2 \,
(V+\alpha)^{-1} \, e^\psi
\end{array}
\end{eqnarray}

\item[(4)]
\begin{eqnarray}
\begin{array}{cl}
& \dps \rho^2 \, (V+\alpha)^{-1} \, \left( |\nabla\psi|^2 +
(\lambda_2-\lambda_1) (1-ce^\psi) \right)\\
\leq & \dps \sqrt{3nc} \, (\lambda_2-\lambda_1) \, \sup \rho^2 \,
(V+\alpha)^{-1} \, \left[ e^\psi \, (1-ce^\psi)
\right]^{\frac{1}{2}}
\end{array}
\end{eqnarray}
\end{enumerate}
\end{theorem}

\vskip 1cm

\section{Estimate of the gap $\lambda_2-\lambda_1$ in terms of the
potential}

\setcounter{equation}{0}

Let $\dps u=\frac{u_2}{u_1}$ be defined as in \S 3. We assume that
it is bounded on $\Omega$ and is zero somewhere in $\Omega$.

Assume that for some $\delta>0$, $u(x_0)= \sup u$ and
$u(x_1)=\delta \sup u$. Then for each smooth function $\rho$ with
compact support and constant $\alpha \geq 0$, we can define
\begin{eqnarray}
L(\rho,\alpha,\delta)=\inf_x \int^1_0 \left( \rho^{-1}
\sqrt{V+\alpha}\right) \left(x(t)\right) |\dot{x}| dt
\end{eqnarray}
where $x$ is any path in $\Omega$ with $x(0)=x_0$ and $x(1)=x_1$.

Now we can define
\begin{eqnarray*}
\begin{array}{rcl}
L_\delta & = & \dps \inf_{\alpha,\rho} L(\rho,\alpha,\delta)
\bigg\{ 20n^2 \sup (V+\alpha)^{-1} \, \rho^2 \, \left(|\Delta\rho|
+
|\nabla\rho|^2\right) \\
& & \hskip 1cm \dps +10n \sup \rho^2 \, (V+\alpha)^{-3} \, \left(
(\Delta V)_+ + |\nabla V|^2 \right) + 10 \frac{\sup
V-\lambda_1}{\sup V+\alpha}\\
& & \hskip 1cm \dps +6n \sup (V+\alpha)^{-\frac{1}{2}} \, \left(
\frac{\sup V-\lambda_1}{\sup V+\alpha} \right)^{\frac{1}{2}} -6n
\inf \left( \rho^2 \, (V+\alpha)^{-1} \, g(x) \right) \bigg\}
\end{array}
\end{eqnarray*}
where $\alpha>0$ is a constant and $\rho$ is any smooth function
with compact support in $\Omega$.

Based on Theorem 3.1, we conclude that
\begin{theorem}
$$\left| \log\frac{1}{\delta} \right| \leq L_\delta +
(\lambda_2-\lambda_1) \left( \frac{1}{\delta} + \inf_{\rho,\alpha}
\left\{ L(\rho,\alpha,\delta) \sup \frac{\rho^2 \,
(V+\alpha)^{-1}}{\delta} \right\} \right)$$
\end{theorem}

In particular if for some $\delta$, $\dps \log
\frac{1-\delta}{\delta} - L_\delta >0$, there is a lower estimate
of $\lambda_2-\lambda_1$, in terms of $L_\delta$ and $\dps
\inf_{\rho,\alpha} \left\{ L(\rho,\alpha,\delta) \sup \left(
\frac{\rho^2 \, (V+\alpha)^{-1}}{\delta} \right) \right\}$.

\vskip 1cm

\section{Oscillation of the function $\dps \frac{u_2}{u_1}$}

\setcounter{equation}{0}

Note that in \S 4, we do not need to assume $u_i$ satisfies any
boundary coordinates on $\Omega$.

If we assume $u_i=0$ on $\ptl\Omega$, $u_1>0$ and
\begin{eqnarray}
\int_\Omega u^2_i=1\\
\int_\Omega u_1\,u_2=0
\end{eqnarray}
we find
\begin{eqnarray}
\int_\Omega u^2\,u^2_1=0\\
\int_\Omega u\,u^2_1=0
\end{eqnarray}

The eigenfunction equations give
\begin{eqnarray}
\int_\Omega |\nabla u_i|^2 + \int_\Omega Vu^2_i=\lambda_i
\end{eqnarray}

Note that assuming (5.1), (5.2) can also be written as
\begin{eqnarray}
\int_\Omega (u_1+u_2)^2=2
\end{eqnarray}
or
\begin{eqnarray}
\int_\Omega (u_1-u_2)^2=2
\end{eqnarray}

Let
\begin{eqnarray}
\Omega_t=\{x\in\Omega | \ d(x_1\,\ptl\Omega) \geq t\}
\end{eqnarray}

We are interested in the behavior of $\dps \int_{\Omega_t} u^2_i$
and $\dps \int_{\Omega_t} u_1\,u_2$. Hence we compute
\begin{eqnarray}
\begin{array}{rcl}
\dps \frac{d^2}{dt^2} \int_{\Omega_t} u^2_i & = & \dps
-\frac{d}{dt} \int_{\delta\Omega_t} u^2_i\\[3mm]
& = & \dps \int_{\ptl\Omega_t} u^2_i \, H_t + 2
\int_{\ptl\Omega_t} u_i \frac{\ptl u_i}{\ptl\nu}\\[3mm]
& = & \dps \int_{\ptl\Omega_t} u^2_i H_t + \int_{\Omega_t} \Delta
(u^2_i)
\end{array}
\end{eqnarray}
where $H_t$ is the mean curvature of $\ptl\Omega_t$, measured by
the outward normal.

\noindent  But
\begin{eqnarray}
\begin{array}{rcl}
\dps \int_{\Omega_t} \Delta(u_i)^2 & = & \dps 2\int_{\Omega_t}
|\nabla
u_i|^2 + 2\int_{\Omega_t} u_i \Delta u_i\\[3mm]
& = & \dps 2 \int_{\Omega_t} (|\nabla u_i|^2 + Vu^2_i)- 2\lambda_i
\int_{\Omega_t} u^2_i\\
& = & \dps 2 \inf (V-\lambda_i)
\end{array}
\end{eqnarray}

\noindent  Since $u_i=0$ on $\ptl\Omega$,
\begin{eqnarray}
\frac{d}{dt} \int_{\Omega_t} u^2_i=0
\end{eqnarray}
where $t=0$.

We conclude that
\begin{eqnarray}
\begin{array}{rcl}
\dps \int_{\Omega_t} u^2_i & = & \dps \int^t_0 \frac{d}{dt} \left(
\int_{\Omega_t} u^2_i \right) + \int_\Omega u^2_i\\[3mm]
& = & \dps \int^t_0 \left( \int^s_0 \frac{d^2}{ds^2} \left(
\int_{\Omega_s} u^2_i \right) \right) + \int_\Omega u^2_i\\[3mm]
& \geq & \dps \inf (V-\lambda_i)t^2 +1
\end{array}
\end{eqnarray}

\noindent Similarly,
\begin{eqnarray}
\begin{array}{rcl}
& & \dps \frac{d^2}{dt^2} \int_{\Omega_t} (u_1+u_2)^2 \\[3mm]
& \geq & \dps 2\int_{\Omega_t} (\nabla u_1+\nabla u_2)^2 +
2\int_{\Omega_t} (u_1+u_2)V(u_1+u_2) -2\int_{\Omega_t}
(u_1+u_2)(\lambda_1\,u_1+\lambda_2\,u_2)\\[3mm]
& \geq & \dps 2\int_{\Omega_t} \left(
V-\frac{\lambda_1+\lambda_2}{2} \right) (u_1+u_2)^2 +
(\lambda_2-\lambda_1) \int_{\Omega_t} u^2_1 +
(\lambda_1-\lambda_2) \int_{\Omega_t} u^2_2\\[3mm]
& \geq & \dps 2\inf \left(V-\frac{(\lambda_1+\lambda_2)}{2}
\right) + (\lambda_1-\lambda_2)
\end{array}
\end{eqnarray}

\begin{eqnarray}
\begin{array}{rcl}
\dps \int_{\Omega_t} (u_1+u_2)^2 & \geq & \dps 2+\left[
\inf_\Omega \left(V-\frac{\lambda_1+\lambda_2}{2}\right) +
\frac{\lambda_1-\lambda_2}{2} \right] t^2\\[3mm]
& = & \dps 2+\left( \inf_\Omega V-\lambda_2 \right) t^2
\end{array}
\end{eqnarray}

\begin{eqnarray}
\begin{array}{rcl}
\dps \int_{\Omega_t} (u_1-u_2)^2 & \geq & \dps 2+\left[
\inf_\Omega \left( V-\frac{\lambda_1+\lambda_2}{2}\right) +
\frac{\lambda_1-\lambda_2}{2} \right] t^2\\[3mm]
& = & \dps 2+ \left(\inf_\Omega V-\lambda_2 \right) t^2
\end{array}
\end{eqnarray}
From (5.14), we obtain
\begin{eqnarray}
\int_{\Omega_t} u_1\,u_2 \geq \frac{1}{2} \left[ \inf_\Omega
\left(V-\lambda_2 \right) \right] t^2
\end{eqnarray}
From (5.15), we obtain
\begin{eqnarray}
\int_{\Omega_t} u_1\,u_2 \leq - \frac{1}{2} \left[ \inf_\Omega
\left(V-\lambda_2 \right) \right] t^2
\end{eqnarray}
Let $\dps u=\frac{u_2}{u_1}$. Then from (5.12),
\begin{eqnarray}
\begin{array}{rcl}
\dps \int_{\Omega_t} u^2\,u^2_1 & = & \dps \int_{\Omega_t}
u^2_2\\[3mm]
& \geq & \dps 1+\inf_\Omega (V-\lambda_2) t^2
\end{array}
\end{eqnarray}
Hence
\begin{eqnarray}
\sup_{\Omega_t} u^2 \geq 1+\inf_\Omega (V-\lambda_2) t^2
\end{eqnarray}
On the other hand,
\begin{eqnarray}
\begin{array}{rcl}
\dps \inf_{\Omega_t} |u| & \leq & \dps \left| \int_{\Omega_t}
u\,u^2_1 \right|\\[3mm]
& \leq & \dps \frac{t^2}{2} \left[ -\inf_\Omega \left( V-\lambda_2
\right) \right]
\end{array}
\end{eqnarray}
Combining (5.19) and (5.20), we conclude
\begin{eqnarray}
\dps \frac{\inf_{\Omega_t}|u|}{\sup_{\Omega_t}|u|} \leq \frac{t^2
\left[-\inf_\Omega \left(V-\lambda_2 \right)
\right]}{2+2\inf_\Omega (V-\lambda_1) t^2}
\end{eqnarray}

In the next section, we shall apply the estimates in \S 4.

\vskip 1cm

\section{Distance function and the estimate of the gap}

\setcounter{equation}{0}

For simplicity, we shall assume that $\Omega$ to be convex in this
section. We also assume that
\begin{eqnarray}
|\nabla V| \leq c_\alpha (V+\alpha)^{\frac{3}{2}}\\
|(\Delta V)_+| \leq c_\alpha (V+\alpha)^3
\end{eqnarray}
We introduce
\begin{eqnarray}
d_\alpha (x_0,x_1)=\inf \int^1_0 \sqrt{V+\alpha} (x(t)) \
|\dot{x}| dt
\end{eqnarray}
where the infinitum is taken over all paths $x: [0,1] \rightarrow
\Omega$ joining $x_0$ to $x_1$.

\noindent If
\begin{eqnarray}
u(x_0)=\sup u
\end{eqnarray}
and
\begin{eqnarray}
u(x_1)=\delta \sup u
\end{eqnarray}

We define $L_\Omega(\delta)$ to be $d(x_0,x_1)$. It is of course
dominated by
\begin{eqnarray}
L(\Omega)=\sup_{\tilde{x_0},\tilde{x_1}\in\Omega}
d(\tilde{x_0},\tilde{x_1})
\end{eqnarray}
Using the terminology of \S 4 and \S 5, we set $\rho$ to be
function of $t$ so that $\rho=0$ on $\ptl\Omega$ and $\rho=1$ when
$\dps t\geq \left(\inf_{\ptl\Omega} V\right)^{-\frac{1}{2}}$.

We can assume
\begin{eqnarray}
\rho^2 \left(|\nabla\log\rho|^2 + |\Delta\log\rho|\right) \leq
3\left(\inf_{\ptl\Omega}V \right)
\end{eqnarray}
Note that
\begin{eqnarray}
\begin{array}{rcl}
\dps \frac{1}{V+\alpha} - \frac{1}{\inf_{\ptl\Omega} V+\alpha} &
\leq & \dps \frac{|\nabla V|}{(V+\alpha)^2}t\\[3mm]
& \leq & \dps \frac{{c_\alpha} t}{(V+\alpha)^{\frac{1}{2}}}\\[3mm]
& \leq & \dps \frac{c_\alpha
t}{(\inf_{\ptl\Omega}V+\alpha)^{\frac{1}{2}}} + \frac{1}{2}
c^2\,t^2
\end{array}
\end{eqnarray}
Hence,
\begin{eqnarray}
\sup (V+\alpha)^{-1} \, \rho^2 (\Delta\rho +|\nabla\rho|^2) \leq
\frac{3(\inf V)}{(\inf_{\ptl\Omega} V+\alpha)} + \frac{c_\alpha
(\inf_{\ptl\Omega} V)^{\frac{1}{2}}}{(\inf_{\ptl\Omega}
V+\alpha)^{\frac{1}{2}}} + \frac{3c^2_\alpha}{2}
\end{eqnarray}
From (5.21), we obtain
\begin{eqnarray}
\frac{\inf_{\Omega_t} |u|}{\sup_{\Omega_t} |u|} \leq
\frac{-\inf_\Omega (V-\lambda_2)}{2\inf_{\ptl\Omega}
V+2\inf_\Omega (V-\lambda_1)}
\end{eqnarray}

Now assume
\begin{eqnarray}
|\inf_\Omega (V-\lambda_2)| \leq \varepsilon\inf_{\ptl\Omega} V
\end{eqnarray}
Then
\begin{eqnarray}
\frac{\inf_{\Omega_t} |u|}{\sup_{\Omega_t} |u|} \leq
\frac{\varepsilon}{2(1-\varepsilon)}
\end{eqnarray}
According to Theorem 4.1, we have proved
\begin{theorem}
Assume (6.1), (6.2) and (6.11) and
$t=(\inf_{\ptl\Omega}V)^{-\frac{1}{2}}$
\begin{eqnarray}
\left| \log \left(\frac{\varepsilon}{2(1-\varepsilon)}\right)
\right| \leq \tilde{c}_\alpha \, L(\Omega_t) +
\frac{2(\lambda_2-\lambda_1)}{\varepsilon} \left[ 1+L(\Omega_t)
\sup_\Omega (V+\alpha)^{-1} \right]
\end{eqnarray}
Here $\tilde{c}_\alpha$ depends on $c_\alpha$ and $-\inf
((V+\alpha)^{-1} \, g(x))$.
\end{theorem}

Note that when $V$ grows fast
$$\left| \log \left(\frac{\varepsilon}{2(1-\varepsilon)} \right)
\right| > \tilde{c}_\alpha \, L(\Omega_+)$$ and we have a lower
bound for $\lambda_2-\lambda_1$ from (6.13).

If we know the location of the points to achieve $\inf_{\Omega_t}
(u)=u(x_0)$ and $\sup_{\Omega_t} |u|=u(x_1)$, then we can replace
$L(\Omega_t)$ by $d_\alpha (x_0,x_1)$. This can be applied when we
have the double well potential.

\end{document}